\documentclass[twoside,12pt]{article}%
\usepackage{amssymb}
\usepackage{amsfonts}
\usepackage{amsmath}
\usepackage{graphicx}%
\providecommand{\U}[1]{\protect\rule{.1in}{.1in}}
\newtheorem{theorem}{Theorem}

\newtheorem{lemma}[theorem]{Lemma}

\newenvironment{proof}[1][Proof]{\noindent\textbf{#1.} }{\ \rule{0.5em}{0.5em}}
\begin{document}

\title{\textbf{Stabilities for Euler-Poisson Equations in Some Special Dimensions }}
\author{Y\textsc{uen} M\textsc{anwai\thanks{E-mail address: nevetsyuen@hotmail.com }}\\\textit{Department of Mathematics and Statistics, Hang Seng School of
Commerce,}\\\textit{Hang Shin Link, Siu Lek Yuen, Shatin,New Territories, Hong Kong}}
\date{Revised 06-Feb-2008}
\maketitle

\begin{abstract}
We study the stabilities and classical solutions of Euler-Poisson equations of
describing the evolution of the gaseous star in astrophysics. In fact, we
extend the study the stabilities of Euler-Poisson equations with or without
frictional damping term to some special dimensional spaces. Besides, by using
the second inertia function in 2 dimension of Euler-Poisson equations, we
prove the non-global existence of classical solutions with $2\int_{\Omega
}(\rho\left\vert u\right\vert ^{2}+2P)dx<gM^{2}-\epsilon$, for any $\gamma$.

\end{abstract}

\section{Introduction}

The evolution of a self-gravitating fluid such as gaseous stars can be
formulated by the Euler-Poisson equations of the following form:
\begin{equation}%
\begin{array}
[c]{rl}%
{\normalsize \rho}_{t}{\normalsize +\nabla\cdot(\rho u)} & {\normalsize =}%
{\normalsize 0,}\\
{\normalsize (\rho u)}_{t}{\normalsize +\nabla\cdot(\rho u\otimes u)+\nabla
P+\beta\rho u} & {\normalsize =}{\normalsize -\rho\nabla\Phi,}\\
S_{t}+u\cdot\nabla S & =0,\\
{\normalsize \Delta\Phi(t,x)} & {\normalsize =\alpha(N)g}{\normalsize \rho,}%
\end{array}
\label{Euler-Poisson}%
\end{equation}
where $\alpha(N)$ is a constant related to the unit ball in $R^{N}$:
$\alpha(1)=2$; $\alpha(2)=2\pi$; For simplicity, we take the constant term
$g=1$. For $N\geq3,$%
\[
\alpha(N)=N(N-2)V(N)=N(N-2)\frac{\pi^{N/2}}{\Gamma(N/2+1)},
\]
where $V(N)$ is the volume of the unit ball in $R^{N}$ and $\Gamma$ is a Gamma
function. As usual, $\rho=\rho(t,x)$, $u=u(t,x)\in\mathbf{R}^{N}$ and $S(t,x)$
are the density, the velocity and the entropy respectively. $P=P(\rho)$\ is
the pressure. In the above system, the self-gravitational potential field
$\Phi=\Phi(t,x)$\ is determined by the density $\rho$ through the Poisson
equation. The equations (\ref{Euler-Poisson})$_{1}$ and (\ref{Euler-Poisson}%
)$_{2}$ are the compressible Euler equation with forcing term. The equation
(\ref{Euler-Poisson})$_{3}$ is the Poisson equation through which the
gravitational potential is determined by the density distribution of the gas
itself. Thus, we called the system (\ref{Euler-Poisson}) the Euler-Poisson
equations. Here, the viscosity term does not appear, that is, the viscous
effect is neglected.\ In this case, the equations can be viewed as a prefect
gas model. For $N=3$, (\ref{Euler-Poisson}) is a classical (nonrelativistic)
description of a galaxy, in astrophysics. See \cite{C}, \cite{M1} for a detail
about the system.

If we take $S(t,x)=\ln K$, for some fixed $K>0$, we have a $\gamma$-law on the
pressure $P(\rho)$, i.e.%
\begin{equation}
{\normalsize P}\left(  \rho\right)  {\normalsize =K\rho}^{\gamma}\doteq
\frac{{\normalsize \rho}^{\gamma}}{\gamma}, \label{gamma}%
\end{equation}
which is a commonly the hypothesis. The constant $\gamma=c_{P}/c_{v}\geq1$,
where $c_{P}$, $c_{v}$\ are the specific heats per unit mass under constant
pressure and constant volume respectively, is the ratio of the specific heats,
that is, the adiabatic exponent in (\ref{gamma}). In particular, the fluid is
called isothermal if $\gamma=1$.

For the physical dimension $N=3$, we are interested in the hydrostatic
equilibrium specified by $u=0,S=\ln K$. According to \cite{C}, the ratio
between the core density $\rho(0)$ and the mean density $\overset{\_}{\rho}$
for $6/5<\gamma<2$\ is given by%
\[
\frac{\overset{\_}{\rho}}{\rho(0)}=\left(  \frac{-3}{z}\overset{\cdot}%
{y}\left(  z\right)  \right)  _{z=z_{0}}%
\]
where $y$\ is the solution of the Lane-Emden equation with $n=1/(\gamma-1)$,%
\[
\overset{\cdot\cdot}{y}+\frac{2}{z}\overset{\cdot}{y}+y^{n}=0,\text{
}y(0)=\alpha>0,\text{ }\overset{\cdot}{y}(0)=0,\text{ }n=\frac{1}{\gamma-1},
\]
and $z_{0}$\ is the first zero of $y(z_{0})=0$. We can solve the Lane-Emden
equation analytically for%
\[
y_{anal}(z)\doteq\left\{
\begin{array}
[c]{ll}%
1-\frac{1}{6}z^{2}, & n=0;\\
\frac{\sin z}{z}, & n=1;\\
\frac{1}{\sqrt{1+z^{2}/3}}, & n=5,
\end{array}
\right.
\]
and for the other values, only numerical values can be obtained. It can be
shown that for $n<5$, the radius of polytropic models is finite; for $n\geq5$,
the radius is infinite.

Gambin \cite{G} and Bezard \cite{B} obtained the existence results about the
explicitly stationary solution $\left(  u=0\right)  $ for $\gamma=6/5:$%
\begin{equation}
\rho=\left(  \frac{3KA^{2}}{2\pi}\right)  ^{5/4}\left(  1+A^{2}r^{2}\right)
^{-5/2}. \label{stationsoluionr=6/5}%
\end{equation}
The Poisson equation (\ref{Euler-Poisson})$_{3}$ can be solved as%
\[
{\normalsize \Phi(t,x)=}\int_{R^{N}}G(x-y)\rho(t,y){\normalsize dy,}%
\]
where $G$ is the Green's function for the Poisson equation in the
$N$-dimensional spaces defined by
\[
G(x)\doteq\left\{
\begin{array}
[c]{ll}%
|x|, & N=1;\\
\log|x|, & N=2;\\
\frac{-1}{|x|^{N-2}}, & N\geq3.
\end{array}
\right.
\]
In the following, we always seek solutions in spherical symmetry. Thus, the
Poisson equation (\ref{Euler-Poisson})$_{3}$ is transformed to%
\[
{\normalsize r^{N-1}\Phi}_{rr}\left(  {\normalsize t,x}\right)  +\left(
N-1\right)  r^{N-2}\Phi_{r}{\normalsize =}\alpha\left(  N\right)
{\normalsize \rho r^{N-1},}%
\]%
\[
\Phi_{r}=\frac{\alpha\left(  N\right)  }{r^{N-1}}\int_{0}^{r}\rho
(t,s)s^{N-1}ds.
\]

In this paper, we concern about existence of global solutions for the
$N$-dimensional Euler-Poisson equations, which may describe the phenomenon
called the core collapsing in physics. And our aim is to construct a family of
such blowup solutions to it. Historically in astrophysics, Goldreich and Weber
\cite{GW} constructed the analytical blowup solution (collapsing) solution of
the $3$-dimensional Euler-Poisson equation for $\gamma=4/3$ for the
non-rotating gas spheres. After that, Makino \cite{M1} obtained the rigorously
mathematical proof of the existence. And in \cite{DXY}, we find the extension
of the above blowup solutions to the case $N\geq3$ and $\gamma=(2N-2)/N$. In
\cite{Y}, the solutions with a form is written as%
\begin{equation}
\left\{
\begin{array}
[c]{c}%
\rho(t,r)=\left\{
\begin{array}
[c]{c}%
\frac{1}{a(t)^{N}}y(\frac{r}{a(t)})^{N/(N-2)},\text{ for }r<a(t)Z_{\mu};\\
0,\text{ for }a(t)Z_{\mu}\leq r.
\end{array}
\right.  \text{, }{\normalsize u(t,r)=}\frac{\overset{\cdot}{a}(t)}%
{a(t)}{\normalsize r,}\text{ }S(t,r)=LnK\\
\overset{\cdot\cdot}{a}(t){\normalsize =}-\frac{\lambda}{a(t)^{N-1}},\text{
}{\normalsize a(0)=a}_{0}>0{\normalsize ,}\text{ }\overset{\cdot}%
{a}(0){\normalsize =a}_{1},\\
\overset{\cdot\cdot}{y}(z){\normalsize +}\frac{N-1}{z}\overset{\cdot}%
{y}(z){\normalsize +}\frac{N(N-2)^{2}V(N)}{(2N-2)K}{\normalsize y(z)}%
^{N/(N-2)}{\normalsize =\mu,}\text{ }y(0)=\alpha>0,\text{ }\overset{\cdot}%
{y}(0)=0,
\end{array}
\right.  \label{solution2}%
\end{equation}
where $\mu=[N(N-2)\lambda]/(2N-2)K$ and the finite $Z_{\mu}$ is the first zero
of $y(z)$.

Recently, in Yuen's results \cite{Y1} there are existing a family of the
blowup solution for the Euler-Poisson equations in radial symmetry in the
$2$-dimensional case,%
\begin{equation}%
\begin{array}
[c]{rl}%
\rho_{t}+u\rho_{r}+\rho u_{r}+{\normalsize \frac{1}{r}\rho u} &
{\normalsize =0,}\\
\rho\left(  u_{t}+uu_{r}\right)  +K\rho_{r} & {\normalsize =-}\frac{2\pi\rho
}{r}\int_{0}^{r}\rho(t,s)sds.
\end{array}
\label{gamma=1}%
\end{equation}
Those are:\bigskip%

\begin{equation}
\left\{
\begin{array}
[c]{c}%
\rho(t,r)=\frac{1}{a(t)^{2}}e^{y(r/a(t))}\text{, }{\normalsize u(t,r)=}%
\frac{\overset{\cdot}{a}(t)}{a(t)}{\normalsize r;}\\
\overset{\cdot\cdot}{a}(t){\normalsize =}-\frac{\lambda}{a(t)},\text{
}{\normalsize a(0)=a}_{0}>0{\normalsize ,}\text{ }\overset{\cdot}%
{a}(0){\normalsize =a}_{1};\\
\overset{\cdot\cdot}{y}(x){\normalsize +}\frac{1}{x}\overset{\cdot}%
{y}(x){\normalsize +\frac{2\pi}{K}e}^{y(x)}{\normalsize =\mu,}\text{
}y(0)=\alpha,\text{ }\overset{\cdot}{y}(0)=0,
\end{array}
\right.  \label{solution1}%
\end{equation}
where $K>0$, $\mu=2\lambda/K$ with a sufficiently small $\lambda$ and $\alpha$
are constants.\newline(1)When $\lambda>0$, the solutions blow up in a finite
time $T$;\newline(2)When $\lambda=0$, if $a_{1}<0$, the solutions blow up at
$t=-a_{0}/a_{1}$.

In this paper, we study the stability of the Euler-Poisson equations with or
without frictional damping in Section 2. After that, we use the second inertia
function in 2 dimension, that is%

\[
{\normalsize H(t)=}\int_{\Omega}{\normalsize \rho(t,x)}\left\vert
{\normalsize x}\right\vert ^{2}{\normalsize dx,}%
\]
to prove the following result:

\begin{theorem}
\label{Thmtamesolution}Let $(\rho(t),u(t))$ be a classical solution with the
property
\[
2\int_{\Omega}(\rho\left\vert u\right\vert ^{2}+2P)dx<gM^{2}-\epsilon,
\]
of the Euler-Poisson equations (\ref{Euler-Poisson}) in 2 dimension on $0\leq
t<T$, here%
\[
M=\int_{\Omega}\rho(t,x)dx,
\]
is the total mass which is constant for any classical solution and $\epsilon$
is an positive uniform constant. Then $T$ must be finite.
\end{theorem}

Initially, in \cite{SI}, Sideris applied the second inertia function to obtain
the instability result for the Euler equations. The similar result for the
Euler equations with frictional damping can be found in \cite{LY}. The
corresponding instability analysis for the $3$-dimensional spherically
symmetric Euler-Poisson Equations for $\gamma\geq4/3$ was appeared in
\cite{MP}. Besides in \cite{DLYY}, the instability similar result of the
$3$-dimensional Euler-Poisson equations for any bounded domain was obtained
later. In a similar way, we extend the results for the case when the some
special dimensions $N\geq2$.

\section{Stability}

In this section, we study the stability of Euler-Poisson equations in higher
dimensions ($N\geq3)$.The total energy can be defined by,%
\begin{equation}
\left\{
\begin{array}
[c]{cc}%
E(t)=\int_{\Omega}\left(  \frac{1}{2}\rho\left\vert u\right\vert ^{2}+\frac
{1}{\gamma-1}P+\frac{1}{2}\rho\Phi\right)  dx, & \text{ for }\gamma>1;\\
E(t)=\int_{\Omega}\left(  \frac{1}{2}\rho\left\vert u\right\vert ^{2}+K\rho
\ln\rho+\frac{1}{2}\rho\Phi\right)  dx, & \text{for }\gamma=1.
\end{array}
\right.  \label{energy1}%
\end{equation}
\newline In \cite{DLYY}, the authors calculated the stationary energy for
$3$-dimensional the Euler-Poisson equations\textbf{.} Here, we generate the
result for $\left(  N\geq2\right)  $-dimensional Euler-Poisson equations.

\begin{lemma}
Suppose the solutions $\rho$ have compact support. For $\gamma>1$ in
Euler-Poisson equations, if $u=0$, the total energy can be expressed by the
integral of the pressure $P$.\newline(1)For $N=2$, we have,%
\[
E(t)=\frac{gM^{2}}{4(\gamma-1)}+\frac{1}{2}\int_{\Omega}\rho(x)\int_{\Omega
}\rho(y)\ln\left\vert x-y\right\vert dydx\text{.}%
\]
(2)For $N\geq3$, we have,%
\[
E(t)=\frac{2N-2-N\gamma}{(\gamma-1)(N-2)}\int_{\Omega}P(x)dx\text{.}%
\]

\end{lemma}

\begin{proof}
Set $u=0$, $\rho(x,t)=\rho(x)$, $\Phi(x,t)=\Phi(x)$.\newline(1)For $N=2$, when
$u=0$ from (\ref{Euler-Poisson})$_{2}$,%
\begin{equation}
{\normalsize \nabla P(x)=-\rho(x)\nabla\Phi(x).} \label{staionary solution}%
\end{equation}
We have, by the divergence theorem%
\begin{align}
2\int P(x)dx  &  =-\int_{\Omega}x\cdot\nabla Pdx\label{Hua}\\
&  =\int_{\Omega}x\cdot\rho(x)\nabla\Phi dx\nonumber\\
&  =g\int_{\Omega}\rho(x)\int_{\Omega}\frac{\rho(y)(x-y)\cdot x}{\left\vert
x-y\right\vert ^{2}}dydx\nonumber\\
&  =gI_{1},\nonumber
\end{align}
where,%
\begin{align*}
I_{1}  &  =\int_{\Omega}\rho(x)\int_{\Omega}\frac{\rho(y)(x-y)\cdot
(x-y)}{\left\vert x-y\right\vert ^{2}}dydx+\int_{\Omega}\rho(x)\int_{\Omega
}\frac{\rho(y)(x-y)\cdot y}{\left\vert x-y\right\vert ^{2}}dydx\\
&  =M^{2}-I_{1}.
\end{align*}
Therefore, we have,%
\[
I_{1}=\frac{1}{2}M^{2},
\]
and%
\[
\int_{\Omega}P(x)dx=\frac{g}{4}M^{2}.
\]
The total energy is%
\begin{align*}
E  &  =\int_{\Omega}\left(  \frac{1}{2}\rho\left\vert u\right\vert ^{2}%
+\frac{1}{\gamma-1}P\right)  dx+\frac{1}{2}\int_{\Omega}\rho\Phi dx\\
&  =\frac{gM^{2}}{4(\gamma-1)}+\frac{1}{2}\int\rho(x)\int\rho(y)\ln\left\vert
x-y\right\vert dydx.
\end{align*}
(2)For $N\geq3$, first we claim the identity,%
\[
N\int_{\Omega}K\rho(x)^{\gamma}dx=\frac{-(N-2)}{2}\int_{\Omega}\rho
(x)\Phi(x)dx.
\]
To prove the claim,%
\[
\rho(x)\nabla_{x}\Phi(x)=-g\rho(x)\int_{\Omega}\nabla_{x}\frac{\rho
(y)}{\left\vert x-y\right\vert ^{N-2}}dy=g\rho(x)\int_{\Omega}\frac
{(N-2)(x-y)\rho(y)}{\left\vert x-y\right\vert ^{N}}dy.
\]
Multiplying this by $x$ and integrating over $\Omega$, we get%
\begin{equation}
\int_{\Omega}x\cdot\rho(x)\nabla\Phi(x)dx=g(N-2)\int\rho(x)\int_{\Omega}%
\frac{\rho(y)(x-y)\cdot x}{\left\vert x-y\right\vert ^{N}}dydx. \label{En11}%
\end{equation}
With the help of divergence theorem, we get%
\[
N\int_{\Omega}K\rho^{r}dx=-\int_{\Omega}x\cdot\nabla\left(  K\rho^{\gamma
}\right)  dx,
\]
and from (\ref{staionary solution}),%
\[
N\int_{\Omega}K\rho^{r}dx=-\int_{\Omega}x\cdot\nabla Pdx=\int_{\Omega}%
x\cdot\rho\nabla\Phi dx.
\]
From (\ref{En11}), we have
\begin{equation}
N\int_{\Omega}K\rho(x)^{r}dx=g(N-2)\int\rho(x)\int_{\Omega}\frac
{\rho(y)(x-y)\cdot x}{\left\vert x-y\right\vert ^{N}}dydx=g\cdot I_{1}\text{,}
\label{En21}%
\end{equation}
where,%
\begin{align*}
I_{1}  &  =(N-2)\int\rho(x)\int_{\Omega}\frac{\rho(y)(x-y)\cdot x}{\left\vert
x-y\right\vert ^{N}}dydx\\
&  =(N-2)\left(  \int\rho(x)\int_{\Omega}\frac{\rho(y)(x-y)\cdot
(x-y)}{\left\vert x-y\right\vert ^{N}}dydx+\int\rho(x)\int_{\Omega}\frac
{\rho(y)(x-y)\cdot y}{\left\vert x-y\right\vert ^{N}}dydx\right) \\
&  =(N-2)\int\rho(x)\int_{\Omega}\frac{\rho(y)}{\left\vert x-y\right\vert
^{N-2}}dydx-I_{1}.
\end{align*}
Thus, we have%
\[
I_{1}=\frac{N-2}{2}\int\rho(x)\int_{\Omega}\frac{\rho(y)}{\left\vert
x-y\right\vert ^{N-2}}dydx.
\]
Hence (\ref{En21}) becomes%
\begin{align*}
N\int_{\Omega}P(x)dx  &  =\frac{g(N-2)}{2}\int_{\Omega}\rho(x)\int_{\Omega
}\frac{\rho(y)}{\left\vert x-y\right\vert ^{N-2}}dydx\\
&  =\frac{N-2}{2}\int_{\Omega}\rho(x)\int_{\Omega}\frac{g\rho(y)}{\left\vert
x-y\right\vert ^{N-2}}dydx\\
&  =-\frac{N-2}{2}\int_{\Omega}\rho(x)\Phi(x)dx.
\end{align*}
The total energy of stationary solutions to Euler-Poisson equations is:%
\begin{align*}
E  &  =\int_{\Omega}\frac{1}{\gamma-1}P(x)dx+\frac{1}{2}\int_{\Omega}%
\rho(x)\Phi(x)dx\\
&  =\int_{\Omega}\frac{1}{\gamma-1}P(x)dx-\frac{N}{N-2}\int_{\Omega}P(x)dx\\
&  =\frac{2N-2-N\gamma}{(\gamma-1)(N-2)}\int_{\Omega}P(x)dx.
\end{align*}
The proof is completed.
\end{proof}

Besides, for the Euler-Poisson equations with frictional damping term, we have
the energy decay estimate lemma:

\begin{lemma}
\label{ChangeRateofEnergy}For the Euler-Poisson equations, suppose the
solutions $(\rho,u)$ have compact support in $\Omega$. With $N\geq3$ and
$\beta\geq0$, we have,
\[
\overset{\cdot}{E}(t)=-\beta\int_{\Omega}\rho\left\vert u\right\vert
^{2}dx\leq0.
\]
Furthermore, if $\beta=0$, then we have,%
\[
\overset{\cdot}{E}(t)=0\text{, i.e. }E(t)\text{ is conserved.}%
\]

\end{lemma}

\begin{proof}
Observe%
\begin{equation}
\left(  \frac{1}{2}\rho\left\vert u\right\vert ^{2}\right)  _{t}=\left(  \rho
u\right)  _{t}\cdot u-\frac{1}{2}\rho_{t}\left\vert u\right\vert
^{2}.\label{energyeq1}%
\end{equation}
Multiplying $-\frac{1}{2}\left\vert u\right\vert ^{2}$ to (\ref{Euler-Poisson}%
)$_{1}$ gives,%
\begin{equation}
-\frac{1}{2}\rho_{t}\left\vert u\right\vert ^{2}-\frac{1}{2}\left\vert
u\right\vert ^{2}\nabla\cdot\left(  \rho u\right)  =0.\label{Energy1}%
\end{equation}
And multiplying $u$ to (\ref{Euler-Poisson})$_{2}$ yields,%
\begin{equation}
\left(  \rho u\right)  _{t}\cdot u+u\cdot\left[  \nabla\cdot\left(  \rho
u\otimes u\right)  \right]  +u\cdot\nabla P=-\rho\nabla\Phi\cdot u-\beta
\rho\left\vert u\right\vert ^{2}.\label{Energy2}%
\end{equation}
We sum up (\ref{Energy1}) and (\ref{Energy2}) to obtain,%
\begin{equation}
\left(  \frac{1}{2}\rho\left\vert u\right\vert ^{2}\right)  _{t}-\frac{1}%
{2}\left\vert u\right\vert ^{2}\nabla\cdot\left(  \rho u\right)
+u\cdot\left[  \nabla\cdot\left(  \rho u\otimes u\right)  \right]
+u\cdot\nabla P=-\rho\nabla\Phi\cdot u-\beta\rho\left\vert u\right\vert
^{2}.\nonumber
\end{equation}
Notice that%
\[
u\cdot\left[  \nabla\cdot\left(  \rho u\otimes u\right)  \right]
=\sum\limits_{i,\text{ }j=1}^{N}u_{i}\left[  \partial_{j}\left(  \rho
u_{j}\right)  u_{i}+\rho u_{j}\partial_{j}u_{i}\right]  =\left\vert
u\right\vert ^{2}\sum\limits_{i=1}^{N}\partial_{i}\left(  \rho u_{i}\right)
+\sum\limits_{i,\text{ }j=1}^{N}\rho u_{i}u_{j}\partial_{j}u_{i}.
\]
Hence, we have,%
\begin{align*}
&  \int_{\Omega}-\frac{1}{2}\left\vert u\right\vert ^{2}\nabla\cdot\left(
\rho u\right)  +u\cdot\left[  \nabla\cdot\left(  \rho u\otimes u\right)
\right]  dx\\
&  =\int_{\Omega}\left[  \frac{1}{2}\left\vert u\right\vert ^{2}%
\sum\limits_{i=1}^{N}\partial_{i}\left(  \rho u_{i}\right)  +\sum
\limits_{i,\text{ }j=1}^{N}\rho u_{i}u_{j}\partial_{j}u_{i}\right]  dx\\
&  =-\int_{\Omega}\left[  \sum\limits_{i=1}^{N}u\cdot\partial_{i}u\rho
u_{i}+\sum\limits_{i,\text{ }j=1}^{N}\rho u_{i}u_{j}\partial_{j}u_{i}\right]
dx\\
&  =-\int_{\Omega}\left[  \sum\limits_{i,\text{ }j=1}^{N}\rho u_{i}%
u_{j}\partial_{i}u_{j}+\sum\limits_{i,\text{ }j=1}^{N}\rho u_{i}u_{j}%
\partial_{j}u_{i}\right]  dx\\
&  =0.
\end{align*}
Besides, the time derivative for the integrand of the second term in
(\ref{energy1})$_{1}$ is:%
\begin{align*}
&  P_{t}=K\gamma\rho^{\gamma-1}\partial_{t}\rho\\
&  =K\gamma\rho^{\gamma-1}\left[  -\sum_{i=1}^{3}\partial_{i}\left(  \rho
u_{i}\right)  \right]  \\
&  =-\sum_{i=1}^{3}K\gamma\rho^{\gamma-1}\partial_{i}\rho u_{i}-\sum_{i=1}%
^{3}K\gamma\rho^{\gamma}\partial_{i}u_{i}\\
&  =-\sum_{i=1}^{3}\partial_{i}Pu_{i}-\sum_{i=1}^{3}\gamma P\partial_{i}%
u_{i}\text{.}%
\end{align*}
Integrating over $\Omega$ and applying integration by part to it, we have%
\[
\int_{\Omega}P_{t}dx=\sum_{i=1}^{3}\left(  \gamma-1\right)  \int_{\Omega
}\partial_{i}Pu_{i}dx.
\]
Thus, we get,
\[
\frac{1}{\gamma-1}\int_{\Omega}P_{t}dx=\int_{\Omega}u\cdot\nabla Pdx.
\]
Applying the divergence theorem and using (\ref{Euler-Poisson})$_{1,3}$, we
have%
\begin{align}
&  -\int_{\Omega}\rho\nabla\Phi\cdot udx\label{Energy4}\\
&  =-\int_{\Omega}\Phi\frac{1}{\alpha(N)g}\Delta\Phi_{t}dx\nonumber\\
&  =-\frac{1}{\alpha(N)g}\int_{\Omega}\Delta\Phi\Phi_{t}dx\nonumber\\
&  =-\int_{\Omega}\rho\Phi_{t}dx.\nonumber
\end{align}
Therefore, from (\ref{Energy4}) and (\ref{Euler-Poisson})$_{1}$, we get%
\[
-\int_{\Omega}\rho\nabla\Phi\cdot udx=-\frac{1}{2}\int_{\Omega}\left(
\Phi\partial_{t}\rho+\rho\partial_{t}\Phi\right)  dx=-\frac{1}{2}\int_{\Omega
}\partial_{t}\left(  \rho\Phi\right)  dx.
\]
By taking the integration for the energy (\ref{energy1}), we have%
\[
\overset{\cdot}{E}(t)=\int_{\Omega}\left(  \frac{1}{2}\rho\left\vert
u\right\vert ^{2}+\frac{1}{\gamma-1}P+\frac{1}{2}\rho\Phi\right)
_{t}dx=-\beta\int_{\Omega}\rho u^{2}dx\leq0.
\]
For $\gamma=1$, we have,%
\[
\int_{\Omega}\rho_{t}dx=-\int_{\Omega}\nabla\left(  \rho u\right)
dx=\int_{\partial\Omega}\rho u\cdot ndx=0,
\]
and%
\begin{align*}
&  \frac{d}{dt}\int_{\Omega}K\rho\ln\rho dx\\
&  =\int_{\Omega}K\rho_{t}\ln\rho dx+\int_{\Omega}K\rho_{t}dx\\
&  =-K\int_{\Omega}\nabla\left(  \rho u\right)  \ln\rho dx\\
&  =\int_{\Omega}u\nabla Pdx.
\end{align*}
Similar calculation gives%
\[
\overset{\cdot}{E}(t)=-\beta\int_{\Omega}\rho\left\vert u\right\vert
^{2}dx\leq0.
\]
The proof is completed.
\end{proof}

In here, we calculate the second derivative of the second inertia function for
the Euler-Poisson equations $(N\geq2)$. 

\begin{lemma}
\label{lem:secondinter}Consider the Euler-Poisson equations without frictional
damping\textbf{,}%
\[%
\begin{array}
[c]{rl}%
{\normalsize \rho}_{t}{\normalsize +\nabla\cdot(\rho u)} & {\normalsize =}%
{\normalsize 0,}\\
{\normalsize (\rho u)}_{t}{\normalsize +\nabla\cdot(\rho u\otimes u)+\nabla P}
& {\normalsize =}{\normalsize -\rho\nabla\Phi,}\\
{\normalsize \Delta\Phi(t,x)} & {\normalsize =}{\normalsize \alpha(N)g\rho,}\\
{\normalsize P} & {\normalsize =}{K{\normalsize \rho}}^{\gamma},
\end{array}
\]
with $(\rho,u)$ is the solution with compact support in $\Omega$, consider the
second inertia, i.e.%
\begin{equation}
{\normalsize H(t)=}\int_{\Omega}{\normalsize \rho(t,x)}\left\vert
{\normalsize x}\right\vert ^{2}{\normalsize dx.}\label{eq21}%
\end{equation}
We have%
\begin{equation}
\left\{
\begin{array}
[c]{cc}%
\overset{\cdot\cdot}{H}(t)=2\int_{\Omega}\left[  (\rho\left\vert u\right\vert
^{2}+NP)dx+\frac{N-2}{2}\rho\Phi\right]  dx, & \text{for }N\geq3;\\
\overset{\cdot\cdot}{H}(t)=2\int_{\Omega}(\rho\left\vert u\right\vert
^{2}+2P)dx-gM^{2}, & \text{for }N=2.
\end{array}
\right.  \label{eq22}%
\end{equation}

\end{lemma}

\begin{proof}
Taking derivatives to (\ref{eq21}) implies:%
\[
\overset{\cdot}{H}(t)=-\int_{\Omega}\nabla\cdot\left(  \rho u\right)
\left\vert x\right\vert ^{2}dx=\int_{\Omega}2x\rho udx,
\]%
\begin{equation}
\overset{\cdot\cdot}{H}(t)=2\int_{\Omega}x\left[  -\nabla\cdot\left(  \rho
u\otimes u\right)  -\nabla P-\rho\nabla\Phi\right]  dx.\label{secondinteria1}%
\end{equation}
We split (\ref{secondinteria1}) into three parts in the following. For the
first part, we apply integration by part,%
\[
-\int_{\Omega}x\cdot\left[  \nabla\cdot\left(  \rho u\otimes u\right)
\right]  dx=-\sum\limits_{j=1}^{N}\int_{\Omega}x_{j}\sum_{i=1}^{N}\partial
_{i}\left(  \rho u_{i}u_{j}\right)  dx=\int_{\Omega}\rho u^{2}dx.
\]
For the second part, we have%
\[
-\int_{\Omega}x\cdot\nabla Pdx=\int_{\Omega}NPdx.
\]
If $N\geq3$, from (\ref{En11}) and (\ref{En21}), we get%
\[
-\int_{\Omega}\rho x\cdot\nabla\Phi dx=\frac{N-2}{2}\int_{\Omega}\Phi\rho dx.
\]
Finally, we have%
\[
\overset{\cdot\cdot}{H}(t)=2\int_{\Omega}\left[  (\rho u^{2}+NP)dx+\frac
{N-2}{2}\rho\Phi\right]  dx.
\]
If $N=2$, since from (\ref{Hua}), we have,%
\[
-\int_{\Omega}\rho x\cdot\nabla\Phi dx=-\frac{1}{2}gM^{2},
\]
it gives%
\[
\overset{\cdot\cdot}{H}(t)=2\int_{\Omega}(\rho\left\vert u\right\vert
^{2}+2P)dx-gM^{2}.
\]
The proof is completed.
\end{proof}

By applying the above lemma, we have the following theorem.

\begin{theorem}
Suppose $(\rho,u)$ is a global classical solution in the Euler-Poisson
equations without frictional damping. We have\newline(1)For $N=3,4$,%
\[
\left\{
\begin{array}
[c]{cc}%
\underset{t\rightarrow\infty}{\lim}\inf\frac{R(t)}{t}\geq\left[  \frac
{(N-2)E}{M}\right]  ^{1/2}, & \text{if }\gamma\geq\frac{2(N-1)}{N}\text{ and
}E>0;\\
\underset{t\rightarrow\infty}{\lim}\inf\frac{R(t)\left\vert \Omega\right\vert
^{(\gamma-1)/2}}{t}{\normalsize \geq}\left[  \frac{(N\gamma-2(N-1))K}%
{\gamma-1}M^{\gamma-1}\right]  ^{1/2}, & \text{if }\gamma>\frac{2(N-1)}%
{N}\text{ and }E=0.
\end{array}
\right.
\]
with $R(t)=\max_{x\in\Omega(t)}\left\{  \left\vert x\right\vert \right\}  $.
Here%
\[
M=\int_{\Omega}\rho(t,x)dx,
\]
is the total mass which is constant for any classical solution.\newline(2)For
$N\geq4$,\newline(2a)the solutions blow up in finite time $T$, if
$1<\gamma\leq2(N-1)/N$ and $E<0$;\newline(2b)the solutions blow up in finite
time $T$, if $1<\gamma<2(N-1)/N$ and $E=0$.
\end{theorem}

\begin{proof}
(1)For $N=3,4$, $\gamma\geq2(N-1)/N$ and $E>0$, we get from Lemma
\ref{lem:secondinter},
\begin{align*}
\overset{\cdot\cdot}{H}(t) &  =2\left\{  \int_{\Omega}\left[  \rho\left\vert
u\right\vert ^{2}+NP\right]  dx+\frac{N-2}{2}\int_{\Omega}\rho\Phi dx\right\}
\\
&  =2\left\{  \int_{\Omega}\left[  \frac{4-N}{2}\rho\left\vert u\right\vert
^{2}+\frac{N\gamma-2(N-1)}{\gamma-1}P\right]  dx+(N-2)E\right\}  \\
&  \geq2(N-2)E>0.
\end{align*}
That is%
\[
{\normalsize H(t)\geq H(0)+}\overset{\cdot}{H}{\normalsize (0)t+(N-2)Et}^{2}.
\]
On the other hand, we have,%
\[
{\normalsize H(0)+}\overset{\cdot}{H}{\normalsize (0)t+(N-2)Et}^{2}%
{\normalsize \leq H}\left(  t\right)  {\normalsize \leq R(t)}^{2}%
{\normalsize M.}%
\]
That is,%
\[
O(\frac{1}{t})+(N-2)E{\normalsize \leq}\frac{R(t)^{2}M}{t^{2}},
\]%
\[
\underset{t\rightarrow\infty}{\lim}\inf\frac{R(t)}{t}\geq\left[  \frac
{(N-2)E}{M}\right]  ^{1/2}.
\]
For $N=3,4$, $\gamma>2(N-1)/N$ and $E=0$, we have
\begin{equation}
{\normalsize M=}\int_{\Omega}{\normalsize \rho dx\leq}\left(  \int_{\Omega
}\rho^{\gamma}dx\right)  ^{1/\gamma}\left\vert \Omega\right\vert
^{(\gamma-1)/\gamma},\label{pressureinequaltiy}%
\end{equation}
and
\begin{align*}
\overset{\cdot\cdot}{H}(t) &  =2\left\{  \int_{\Omega}\left[  \frac{4-N}%
{2}\rho\left\vert u\right\vert ^{2}+\frac{N\gamma-2(N-1)}{\gamma-1}P\right]
dx+(N-2)E\right\}  \\
&  \geq2\int_{\Omega}\frac{N\gamma-2(N-1)}{\gamma-1}Pdx.
\end{align*}
Thus, we get,%
\[
\overset{\cdot\cdot}{H}(t){\normalsize \geq}\frac{2\left[  N\gamma
-2(N-1)\right]  K}{\gamma-1}\left\vert \Omega\right\vert ^{1-\gamma
}{\normalsize M}^{\gamma}{\normalsize >0,}%
\]%
\[
{\normalsize H(0)+}\overset{\cdot}{H}(0){\normalsize t+}\frac{\left[
N\gamma-2(N-1)\right]  K}{\gamma-1}\left\vert \Omega\right\vert ^{1-\gamma
}{\normalsize M}^{\gamma}{\normalsize t}^{2}{\normalsize \leq H(t)\leq
R(t)}^{2}{\normalsize M,}%
\]%
\[
{\normalsize O(}\frac{1}{t}{\normalsize )+}\frac{\left[  N\gamma
-2(N-1)\right]  K}{\gamma-1}\left\vert \Omega\right\vert ^{1-\gamma
}{\normalsize M}^{\gamma}{\normalsize \leq}\frac{R(t)^{2}M}{t^{2}}.
\]
This gives%
\[
\underset{t\rightarrow\infty}{\lim}\inf\frac{R(t)\left\vert \Omega\right\vert
^{(\gamma-1)/2}}{t}{\normalsize \geq}\left[  \frac{(N\gamma-2(N-1))K}%
{\gamma-1}M^{\gamma-1}\right]  ^{1/2}.
\]
(2a)For $N\geq4$, $\gamma\leq2(N-1)/N$ and $E<0$, we have,%
\begin{align*}
\overset{\cdot\cdot}{H}(t) &  =2\left\{  (N-2)E+\int_{\Omega}\left[
\frac{4-N}{2}\rho\left\vert u\right\vert ^{2}+\frac{N\gamma-2(N-1)}{\gamma
-1}P\right]  dx\right\}  \\
&  \leq2(N-2)E.
\end{align*}
That is,%
\[
H(t)\leq H(0)+\overset{\cdot}{H}(0)t+(N-2)Et^{2}.
\]
As $E<0$, $H(0)+\overset{\cdot}{H}(0)t+(N-2)Et^{2}$ will become negative if
$t$ is sufficiently large. But, $H(t)$ is always non-negative. Thus we have a
contradiction. When $N\geq4$, $\gamma\leq2(N-1)/N$ and $E<0$ the classical
solutions blow up after some finite time $T$.

(2b)For $N\geq4$, $\gamma<2(N-1)/N$ and $E=0$, we obtain,%
\begin{align*}
\overset{\cdot\cdot}{H}(t) &  =2\left\{  (N-2)E+\int_{\Omega}\left[
\frac{4-N}{2}\rho\left\vert u\right\vert ^{2}+\frac{N\gamma-2(N-1)}{\gamma
-1}P\right]  dx\right\}  \\
&  \leq2\int_{\Omega}\frac{N\gamma-2(N-1)}{\gamma-1}Pdx.
\end{align*}
Form (\ref{pressureinequaltiy}), we get,%
\[
\overset{\cdot\cdot}{H}(t){\normalsize \leq}\frac{2(N\gamma-2(N-1))K}%
{\gamma-1}\left\vert \Omega\right\vert ^{1-\gamma}{\normalsize M}^{\gamma}.
\]
After solving it, we have,%
\[
H(t)\leq H(0)+\overset{\cdot}{H}(0)t+\frac{(N\gamma-2(N-1))K}{\gamma
-1}\left\vert \Omega\right\vert ^{1-\gamma}M^{\gamma}t^{2}.
\]
As $\gamma<2(N-1)/N$, we use the above same argument in (3a) to prove
non-global existence of the solution exists. Thus, when $N\geq4$,
$\gamma<2(N-1)/N$ and $E=0$ the classical solutions blow up after some finite
time $T$. The proof is completed.
\end{proof}

By applying the similar method in the above section to the Euler-Poisson
equations with frictional damping $\beta>0$,
\begin{equation}%
\begin{array}
[c]{rl}%
{\normalsize \rho}_{t}{\normalsize +\nabla\cdot(\rho u)} & {\normalsize =}%
{\normalsize 0,}\\
{\normalsize (\rho u)}_{t}{\normalsize +\nabla\cdot(\rho u\otimes u)+\nabla
P+} & {\normalsize =}{\normalsize -\rho\nabla\Phi-\beta\rho u,}\\
{\normalsize \Delta\Phi(t,x)} & {\normalsize =}{\normalsize \alpha(N)g\rho,}\\
{\normalsize P} & {\normalsize =}{K{\normalsize \rho}}^{\gamma}.
\end{array}
\label{DD}%
\end{equation}

\begin{theorem}
Suppose $(\rho,u)$\ is a global solution of the Euler-Poisson equations with
frictional damping,\textbf{ }i.e. (\ref{DD}), for $N\geq4$,\newline(1)the
solutions blow up in finite time $T$, if $1<\gamma\leq2(N-1)/N$ and
$E(0)<0$;\newline(2)the solutions blow up in finite time $T$, if
$1<\gamma<2(N-1)/N$ and $E(0)\leq0$.
\end{theorem}

\begin{proof}
(1)For $N\geq4$, $\gamma\leq2(N-1)/N$ and $E(0)<0$, we get,%
\begin{align*}
\overset{\cdot}{H}(t)  &  =-\frac{1}{\beta}\overset{\cdot\cdot}{H}(t)+\frac
{2}{\beta}\left\{  \int_{\Omega}\left[  \frac{4-N}{2}\rho\left\vert
u\right\vert ^{2}+\frac{N\gamma-2(N-1)}{\gamma-1}P\right]  dx+(N-2)E\right\}
\\
&  \leq-\frac{1}{\beta}\overset{\cdot\cdot}{H}(t)+\frac{2(N-2)}{\beta}E(t).
\end{align*}
From Lemma \ref{ChangeRateofEnergy}, we have,%
\[
\overset{\cdot}{H}(t)\leq-\frac{1}{\beta}\overset{\cdot\cdot}{H}%
(t)+\frac{2(N-2)E(0)}{\beta}.
\]
That is%
\[
H(t)\leq C_{1}+C_{2}e^{-\beta t}+\frac{2(N-2)E(0)}{\beta}t.
\]
As $E(0)<0$, $C_{1}+C_{2}e^{-\beta t}+\left[  2(N-2)E(0)t\right]  /\beta$ will
become negative if $t$ is sufficiently large. But, $H(t)$ is always
non-negative. Thus we have a contradiction. When $N\geq4$, $\gamma
\leq2(N-1)/N$ and $E(t)<0$ for $t\geq0$, the classical solutions blow up after
some finite time $T$.

(2)For $N\geq4$, $\gamma<2(N-1)/N$ and $E(0)=0$, we get%
\begin{align*}
\overset{\cdot}{H}(t)  &  =-\frac{1}{\beta}\overset{\cdot\cdot}{H}(t)+\frac
{2}{\beta}\left\{  \int_{\Omega}\left[  \frac{4-N}{2}\rho\left\vert
u\right\vert ^{2}+\frac{N\gamma-2(N-1)}{\gamma-1}P\right]  dx+(N-2)E\right\}
\\
&  \leq-\frac{1}{\beta}\overset{\cdot\cdot}{H}(t)+\frac{2}{\beta}\int_{\Omega
}\frac{N\gamma-2(N-1)}{\gamma-1}Pdx.
\end{align*}
This implies%
\[
H(t)\leq C_{1}+C_{2}e^{-\beta t}+\left[  \frac{2}{\beta}\int_{\Omega}%
\frac{N\gamma-2(N-1)}{\gamma-1}Pdx\right]  t.
\]
Since%
\[
{\normalsize M=}\int_{\Omega}{\normalsize \rho dx\leq}\left(  \int_{\Omega
}\rho^{\gamma}dx\right)  ^{1/\gamma}\left\vert \Omega\right\vert
^{(\gamma-1)/\gamma},
\]
we get,%
\[
H(t)\geqslant C_{1}+C_{2}e^{-\beta t}+\left[  \frac{2(N\gamma-2(N-1))K}%
{\beta(\gamma-1)}M^{\gamma}\left\vert \Omega\right\vert ^{1-\gamma}\right]
t.
\]
As $\gamma<2(N-1)/N$, we use the same argument in (1) to give the non-global
existence of the classical solution. Thus, when $N\geq4$, $\gamma<2(N-1)/N$
and $E(t)\leq0$ for $t\geq0$, the classical solutions blow up after some
finite time $T$. The proof is completed.
\end{proof}

By the similar argument, when $\gamma=\frac{2(N-1)}{N}$ for $N\geq4$, any
perturbation with negative total energy leads to non-global existence.

\section{Non-global Existence for 2-d Classical Solutions}

We present the proof of the theorem \ref{Thmtamesolution} here.

\begin{proof}
From the second inertia function in 2 dimension for the Euler-Poisson
equations (\ref{eq22})$_{2}$ without frictional damping $(\beta=0)$, we have%
\[
\overset{\cdot\cdot}{H}(t)=2\int_{\Omega}(\rho\left\vert u\right\vert
^{2}+2P)dx-gM^{2}.
\]
Under the assumption in Theorem \ref{Thmtamesolution}, we get%
\[
\overset{\cdot\cdot}{H}(t)<-\epsilon,
\]%
\[
H(t)=-\frac{\epsilon}{2}t^{2}+d_{0}t+d_{1},
\]
where $d_{0}$ and $d_{1\text{ }}$are some constants.

After some finite time, we meet a contradiction%
\[
H(t)<0,
\]

as $H(t)$ cannot be negative. Thus the classical solutions of such kinds in
global existence are not possible. The proof is completed.
\end{proof}


\begin{thebibliography}{99}                                                                                               %


\bibitem {B}M. Bezard, \textit{Existence locale de solutions pour les
equations d'Euler-Poisson. (French) [Local Existence of Solutions for
Euler-Poisson Equations]} Japan J. Indust. Appl. Math. \textbf{10} (1993), no.
3, 431--450.

\bibitem {C}S. Chandrasekhar, \textit{An Introduction to the Study of Stellar
Structure}, Univ. of Chicago Press, 1939.

\bibitem {DLYY}Y.B. Deng, T.P. Liu, T. Yang and Z.A. Yao, \textit{Solutions of
Euler-Poisson Equations for Gaseous Stars}, Arch. Rational Mech. Anal.
\textbf{164 }(2002), 261-285.

\bibitem {DXY}Y.B. Deng, J.L. Xiang and T. Yang, \textit{Blowup Phenomena of
Solutions to Euler-Poisson Equations}, J. Math. Anal. Appl. \textbf{286}
(1)(2003), 295-306.

\bibitem {G}P. Gamblin, \textit{Solution reguliere a temps petit pour
l'equation d'Euler-Poisson. (French) [Small-time Regular Solution for the
Euler-Poisson Equation]} Comm. Partial Differential Equations \textbf{18}
(1993), no. 5-6, 731--745.

\bibitem {GW}P.Goldreich, S. Weber, \textit{Homologously Collapsing Stellar
Cores}, Astrophys, J. \textbf{238}, 991 (1980).

\bibitem {LY}T.P. Liu and T. Yang, \textit{Compressible Euler Equations with
Vacuum}, J. Differential Equations \textbf{140} (1997), No. 2, 223-237.

\bibitem {M1}T. Makino, \textit{Blowing up Solutions of the Euler-Poission
Equation for the Evolution of the Gaseous Stars}, Transport Theory and
Statistical Physics \textbf{21} (1992), 615-624.

\bibitem {MP}T. Makino and B. Perthame, \textit{Sur les Solutions a symmetric
spherique de lequation d'Euler-poisson Pour levolution d'etoiles
gazeuses,(French) [On Radially Symmetric Solutions of the Euler-Poisson
Equation for the Evolution of Gaseous Stars], }Japan J. Appl. Math. \textbf{7}
(1990), 165-170.

\bibitem {SI}T.C. Sideris, \textit{Formation of Singularities in
Three-dimensional Compressible Fluids}, Comm. Math. Phys. \textbf{101} (1985),
No. 4, 475-485.

\bibitem {Y}M.W. Yuen, \textit{Blowup Solutions for a Class of Fluid Dynamical
Equations in }$R^{N}$ , J. Math. Anal. Appl. \textbf{329} (2)(2007), 1064-1079.

\bibitem {Y1}M.W. Yuen, Analytical Blowup Solutions to the 2-dimensional
Isothermal Euler-Poisson Equations of Gaseous Stars, J. Math. Anal. Appl.
\textbf{341 (}1\textbf{)(}2008\textbf{), }445-456.
\end{thebibliography}
\end{document}